# Locally adaptive estimation methods with application to univariate time series[*]

Mstislav Elagin

*Weierstrass Institute for Applied Analysis and Stochastics*
*Mohrenstr. 39, 10117 Berlin, Germany*
*e-mail:* elagin@wias-berlin.de

**Abstract:** The paper offers a unified approach to the study of three locally adaptive estimation methods in the context of univariate time series from both theoretical and empirical points of view. A general procedure for the computation of critical values is given. The underlying model encompasses all distributions from the exponential family providing for great flexibility. The procedures are applied to simulated and real financial data distributed according to the Gaussian, volatility, Poisson, exponential and Bernoulli models. Numerical results exhibit a very reasonable performance of the methods.

**AMS 2000 subject classifications:** Primary 62M10; secondary 62F10, 62P20.
**Keywords and phrases:** Adaptive estimation, Local homogeneity, Model selection, Stagewise aggregation, Volatility model, Poisson model, Exponential model, Bernoulli model, Propagation, Oracle.

## Introduction

Since ARMA (Box and Jenkins, 1994), the classical linear model for stationary time series, a multitude of extensions for the nonlinear and non-stationary case has been developed. Models for nonlinear time series can be divided into two groups: those modelling nonlinearity in the conditional mean, e.g. threshold auto-regressive model (Tong and Lim, 1980), and those modelling nonlinearity in the variance, e.g. conditional heteroskedastic (CH) models including ARCH (Engle, 1982) and GARCH (Bollerslev, 1986) and their extensions that have been inspired by the particular behaviour of financial time series ("stylized facts") such as volatility clustering, heavy tails, asymmetry etc. Generalizations for the non-stationary case include ARIMA and ARFIMA as well as IGARCH (Nelson, 1990) and FIGARCH Baillie et al. (1996) models for (fractionally) integrated time series. Fan and Wenyang (2008) give an overview on the modern developments in the area of varying coefficient models originated by Hastie and Tibshirani (1993).

However, the question of whether ubiquitous financial time series are nonstationary, e.g. possess a unit root, or experience structural breaks has been ex-tensively debated (Nelson and Plosser, 1982), (Perron, 1989). It has been shown (Mikosch and Stărică, 2004) that long range memory effects in financial time series may be caused by structural breaks rather than result from nonstationarity. Further it has been argued (Diebold and Inoue, 2001) that structural breaks can easily be overlooked with negative impact on the quality of modelling, estimation and forecasting (Hillebrand, 2005). Perron (1989) introduced a model with breaks that are exogenous, i.e. induced by external influence at known times, Zivot and Andrews (2002) suggested a sequential method of testing for endogenous breaks, i.e. not known beforehand. In this paper we also test for endogenous breaks sequentially and apply locally parametric estimation methods to stationary subsets of the series. Methods of this kind have been presented e.g. in Fan and Gu (2003) for adaptive selection of the decay factor used to weight components of the pseudo-likelihood function, in Dahlhaus and Subba Rao (2006) for the formulation of the locally stationary ARCH($\infty$) processes, in Cheng et al. (2003) for locally choosing parameters of a filter.

---

[*]Work supported by the Collaborative Research Centre 649 "Economic Risk" funded by the Deutsche Forschungsgemeinschaft.

The mentioned locally parametric methods are the local change point (LCP) procedure (Mercurio and Spokoiny, 2004), the local model selection (LMS), also known as the intersection of confidence intervals (ICI) (Katkovnik and Spokoiny, 2008), and the stagewise aggregation (SA) (Belomestny and Spokoiny, 2007). Until present the properties of these procedures were proven for the corresponding particular cases. In this paper we present a unified approach to adaptive estimation upon which the three aforementioned methods are based, treat theoretical properties of the estimators in a general way way and give a universal procedure for the choice of parameters (*critical values*).

Another focus of the paper is on the generality of the model. Modelling financial data involves the use of various distributions. For instance, the volatility model is applied to study the dynamics of squared returns, transactions on the market can be looked at as Bernoulli random events, time between transactions is usually assumed to follow the exponential or Weibull distribution, and transaction intensity can be described by the Poisson distribution. In this paper a rather general model allowing for an arbitrary distribution from the exponential family is adopted. Models of this kind, however with accent of generalized linear models, have been used to study binary, categorical and count time series (Fokianos and Kedem, 2003).

The paper is organized as follows. Section 1 is devoted to the formulation of the problem in the context of univariate time series and theoretical introduction. Section 2 describes the local parametric methods. In the Section 3 the procedure to obtain *critical values*, essential parameters of the procedures is given. Section 4 collects theoretical properties of the procedures, including the main result, the oracle inequality, stating that the quality of the adaptive estimates provided by the procedures is comparable with the quality of the best possible estimate. Finally, Section 5 demonstrates the performance of the procedures on simulated and real data belonging to various distributions.

## 1. Model and set-up

Consider a one-dimensional stochastic process $Y_t$ in discrete time $t \in \mathbb{N}$ that is progressively measurable w.r.t. a filtration $(\mathcal{F}_t)$. We assume that the process $Y_t$ conditionally on the past has a distribution $P_v$ from a given family of distributions $\mathcal{P}$ indexed by a parameter $v$:

$$\mathcal{L}(Y_t|\mathcal{F}_{t-1}) \in \mathcal{P} = (P_v), \quad v \in \Upsilon.$$

In the general case, the parameter $v$ is a predictable stochastic process, leading to the following setup:
$$Y_t|\mathcal{F}_{t-1} \sim P_{v_t}. \tag{1}$$

Our goal is to infer on $v_t$.

### 1.1. Parametric and local parametric estimation and inference

In the parametric case the parameter $v_t$ is assumed to depend on time through a time-varying function $f_t(\theta)$ of an argument $\theta$:

$$v_t = f_t(\theta). \tag{2}$$

Then the estimation of $v_t$ based on the observations $Y_t$ can be carried out using the maximum likelihood method. The model (1)–(2) leads to the log-likelihood function

$$L(\theta) = \sum_t \ell(Y_t, v_t) = \sum_t \ell(Y_t, f_t(\theta))$$

where $\ell$ is the log-density of the distribution family $\mathcal{P}$. The estimate of the unknown parameter is then obtained by maximizing the log-likelihood function w.r.t. $\theta$:

$$\tilde{v}_t = f_t(\tilde{\theta}) \quad \text{with} \quad \tilde{\theta} = \arg\max_\theta L(\theta) = \arg\max_\theta \sum_t \ell(Y_t, f_t(\theta)).$$

However, the parametric assumption is often too restrictive in the practice. We therefore follow the *local parametric approach*: it is supposed that for the time point of estimation $t^\diamond$ there exists some interval $I = [t^\diamond - N_I, t^\diamond]$ of length $N_I$ within which the parametric regression function $f_t(\theta)$ describes the process adequately. We aim to determine that interval.

### 1.2. Local constant approach for the exponential family

Below we consider a special case of the parametric regression function, namely the *local constant parametric* one:

$$f_t(\theta) \equiv \theta, \, t \in I,$$

where both the length $N_I$ of the interval and the value $\theta$ are to be estimated from the data. The likelihood function then assumes the form

$$L_I(\theta) = \sum_{t \in I} \ell(Y_t, \theta).$$

Further, we limit the discussion to a particular (*exponential*) family of distributions allowing for simple expression of the likelihood function and of the maximum likelihood estimate. Recall that distributions belonging to the exponential family possess the densities of the form

$$p(y, \theta) = p(y) e^{y \upsilon(\theta) - B(\theta)}, \, \theta \in \Theta, \, y \in \mathcal{Y}, \tag{3}$$

where $\upsilon(\theta)$, $B(\theta)$ and and $p(y)$ are some given functions.

One distinguishes between the *natural* and *canonical* parametrisations. A parameter $\theta$ is called natural if it satisfies the equality

$$\mathbf{E}_\theta Y = \int y p(y, \theta) P(dy) = \theta$$

for all $\theta \in \Theta$. Hence, the local maximum likelihood estimate of the parameter is the average of the observations:

$$\tilde{\theta}_I = \arg\max_{\theta \in \Theta} L_I(\theta) = \sum_{t \in I} Y_t / N_I.$$

A parameter $v$ is canonical if the density $p$ of any measure $P_v \in \mathcal{P}$ w.r.t. the dominating measure $P$ can be represented as

$$p(y, v) \equiv \frac{dP_v}{dP}(y) = p(y) e^{yv - d(v)},$$

where $d(v)$ is some convex function on $\Upsilon$. Canonical parametrisation is convenient due to the parameter entering the log likelihood linearly and allowing for mixing of distributions. This property is applied in the method of stagewise aggregation (see Section 2.3).

For the exposition we will need a measure of difference between distributions $P_{\theta_1}$ and $P_{\theta_2}$ from the family $\mathcal{P}$ indexed by $\theta_1$ and $\theta_2$. One popular measure is the Kullback – Leibler divergence defined as follows:

$$\mathcal{K}(\theta_1, \theta_2) = \mathbf{E}_{\theta_1} \log \frac{dP_{\theta_1}}{dP_{\theta_2}}.$$

Its form for the distributions used in this paper is given in the Table 1. Although the form of the Kullback – Leibler divergence depends on the parametrisation, it is true that

$$\mathcal{K}_\theta(\theta_1, \theta_2) = \mathcal{K}_\upsilon\left(\upsilon(\theta_1), \upsilon(\theta_2)\right),$$

for the natural parametrisation, $\mathcal{K}_\upsilon(\cdot, \cdot)$ is the form for the canonical parametrisation, and $\upsilon(\theta)$ is the function giving the relation between canonical and natural parameters. Due to this property, many results of this section hold for both parametrisations.

| Model | $\mathcal{K}(\theta, \theta')$ |
|---|---|
| Gaussian regression | $(\theta - \theta')^2/(2\sigma^2)$ |
| Bernoulli | $\theta \log(\theta/\theta') + (1-\theta)\log\{(1-\theta)/(1-\theta')\}$ |
| Poisson model | $\theta \log(\theta/\theta') - (\theta - \theta')$ |
| Exponential model | $\theta/\theta' - 1 - \log(\theta/\theta')$ |
| Volatility model | $^1/_2(\theta/\theta' - 1) - ^1/_2 \log(\theta/\theta')$ |

TABLE 1
*Kullback – Leibler divergence for some distributions from the exponential family with natural parametrisation.*

As the following Theorem 1 asserts, the *fitted likelihood* $L(\tilde{\theta}, \theta)$, i.e. the log ratio of the maximum value of the likelihood function to its value at an arbitrary point $\theta$, is closely related to the Kullback – Leibler divergence.

**Theorem 1 (Polzehl and Spokoiny (2006)).** *For distributions from the exponential family it is true that*

$$L_I(\tilde{\theta}_I, \theta) = \max_{\theta'} L_I(\theta', \theta) = L_I(\tilde{\theta}_I) - L_I(\theta) = N_I \mathcal{K}(\tilde{\theta}_I, \theta). \tag{4}$$

The following Theorem 2 claims that in the parametric case the estimation loss measured by $L_I(\tilde{\theta}_I, \theta^*)$ is with high probability bounded.

**Theorem 2 (Polzehl and Spokoiny (2006), Theorem 2.1).** *Let the observations $Y_t$ be i.i.d. from $P_{\theta^*}$, $\theta^* \in \Theta$. Then for any $\mathfrak{z} > 0$*

$$\mathbf{P}_{\theta^*}\left\{L_I(\tilde{\theta}_I, \theta^*) > \mathfrak{z}\right\} \leq 2e^{-\mathfrak{z}},$$

*where $\mathbf{P}_{\theta^*}$ is the joint distribution.*

Based on this result, one can bound the risk (expected loss) for a power loss function with the power $r$:

**Theorem 3.** *Let the observations $Y_t$ be i.i.d. according to $P_{\theta^*}$. Denote as $\mathcal{R}_{r,\theta^*,I}$ the risk of estimation based on the observations from the interval $I$:*

$$\mathcal{R}_{r,\theta^*,I} = \mathbf{E}_{\theta^*} \left| L_I(\tilde{\theta}_I, \theta^*) \right|^r \tag{5}$$

*with some $r > 0$. Then the estimation risk is bounded:*

$$\mathcal{R}_{r,\theta^*,I} \leq \mathfrak{r}_r, \quad where \quad \mathfrak{r}_r = 2r \int_{\mathfrak{z} \geq 0} \mathfrak{z}^{r-1} e^{-\mathfrak{z}} d\mathfrak{z} = 2r\Gamma(r).$$

**Proof.** Denote the event $\{L_I(\tilde{\theta}_I, \theta^*) > \mathfrak{z}\}$ by $\mathbb{A}$. Using the expectation representation formula, integrating by parts and applying Theorem 2 yields

$$\mathcal{R}_{r,\theta^*,I} = -\int_{\mathfrak{z} \geq 0} \mathfrak{z}^r d\mathbf{P}_{\theta^*}\{\mathbb{A}\} = r \int_{\mathfrak{z} \geq 0} \mathfrak{z}^{r-1} \mathbf{P}_{\theta^*}\{\mathbb{A}\} d\mathfrak{z} \leq 2r \int_{\mathfrak{z} \geq 0} \mathfrak{z}^{r-1} e^{-\mathfrak{z}} d\mathfrak{z}. \quad \square$$

It follows from the Theorem 2 that if some $\mathfrak{z}_\alpha$ is such that $2e^{-\mathfrak{z}_\alpha} \leq \alpha$, then

$$\mathcal{E}_I(\mathfrak{z}_\alpha) = \left\{\theta : L_I(\tilde{\theta}_I, \theta) \leq \mathfrak{z}_\alpha\right\}$$

is the $\alpha$-confidence set for the parameter $\theta$.

### 1.3. Nearly parametric case

In practice the parametric assumption may be overly stringent and not hold even within an arbitrarily small interval. However, it turns out that results of Theorem 3 still hold to a certain extent even if the parametric assumption is violated. Deviation from the parametric case can be described by a magnitude

$$\Delta_I(\theta) = \sum_{t \in I} \mathcal{K}(v_t, \theta)$$

which is in general random. The "nearly parametric" case occurs whenever the following condition is fulfilled:

**Condition 1 (Small modelling bias).** *For a given interval $I^\circ$ there exists a parameter value $\theta^\circ \in \Theta$ such that the expectation under the true measure of the divergence $\Delta_{I^\circ}(\theta^\circ)$ is bounded by some $\Delta^\circ \geq 0$:*

$$\mathbf{E}\Delta_{I^\circ}(\theta^\circ) \leq \Delta^\circ.$$

The following theorem generalizes the result of the Theorem 3 to the nearly parametric case:

**Theorem 4.** *Let the small modelling bias condition hold. Then for $r > 0$*

$$\mathbf{E}\log\left(1 + \frac{\left|L_{I^\circ}(\tilde{\theta}_{I^\circ}, \theta^\circ)\right|^r}{\mathcal{R}_{r,\theta^\circ,I^\circ}}\right) \leq 1 + \Delta^\circ,$$

*where $\mathcal{R}_{r,\theta^\circ,I^\circ}$ is defined analogously to* (5).

**Proof.** It suffices to apply Lemma 2 (information-theoretical bound) with

$$\zeta = \left|L_{I^\circ}(\tilde{\theta}_{I^\circ}, \theta^\circ)\right|^r / \mathcal{R}_{r,\theta^\circ,I^\circ},$$

$\mathbf{P} = \mathbf{P}_{I^\circ}$ (true measure on the interval $I^\circ$), $\mathbf{P}_0 = \mathbf{P}_{\theta^\circ}$ (parametric measure with parameter $\theta^\circ$) and employ that $\mathbf{E}_0\zeta = 1$ by definition of $\mathcal{R}_{r,\theta^\circ,I^\circ}$. □

This result means that the risk of the maximum likelihood estimate in the nearly parametric case is comparable to that in the parametric case. The logarithm under the expectation is induced by the information-theoretical bound, and $\Delta^\circ$ can be interpreted as payment for violation of the parametric assumption. Theorem 4 leads to the notion of the *oracle* estimate as the one corresponding to the largest interval under the small modelling bias condition. Performance of the oracle estimate will be the basis for measuring the performance of adaptive estimates in the following sections.

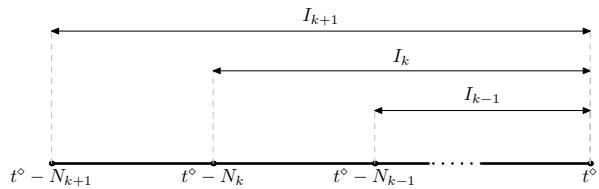

Fig 1. *Nested intervals.*

## 2. Adaptive estimation methods

In the Theorem 4 the oracle estimate has been defined as the estimate corresponding to the largest interval under the small modelling bias condition. Unfortunately, since the oracle depends on the unknown true measure, it is not available and can only be mimicked. Based on a sequence of growing intervals and associated estimates, we aim to construct an adaptive estimate matching the oracle quality. Let $t^\diamond$ denote the time point at which the parameter $\upsilon$ is to be estimated, and $\{I_k\}_{k=1}^{K}$ be a sequence of growing intervals of length $N_k$ with the common right edge at $t^\diamond$ (Fig. 1): $I_k = [t^\diamond - N_k, t^\diamond[$. We associate with each interval $I_k$ from this sequence the corresponding maximum likelihood estimate, which we shall call *weak estimate*. To simplify the notation, in the sequel we denote $\tilde{\theta}_k = \tilde{\theta}_{I_k}$, $L_k(\theta_1, \theta_2) = L_{I_k}(\theta_1, \theta_2)$ etc. The *model selection*-type approach to adaptive estimation consists in picking the estimate corresponding to the largest interval satisfying a certain predicate $\mathcal{A}$ together with all of the smaller intervals:

$$\hat{\theta} = \tilde{\theta}_{\hat{k}} \text{ with } \hat{k} = \max\{k : \mathcal{A}(k) = 1 \text{ and } \mathcal{A}(k-1) = 1\} \qquad (6)$$

where $\mathcal{A}(1) = 1$ by definition. Below we present two model selection-type procedures specifying the predicate $\mathcal{A}$.

In contrast, the *aggregation* approach gradually builds up the adaptive estimate by taking convex combinations of the weak estimates. It will be shown that adaptive estimates obtained by all three procedures exhibit the risk comparable to the oracle risk.

Note that all procedures are pointwise, i.e. they provide an estimate $\hat{\theta} = \hat{\theta}(t)$ for every time point $t$.

### 2.1. Local change point detection

The method of local change point (LCP) detection introduced in Mercurio and Spokoiny (2004) is a procedure for discovering a change point within an interval provided there can be at most one such point. Suppose that an interval $\mathcal{I}'' = (\tau, t^\diamond]$ has been found to contain no change points (Fig. 2). To test the null hypothesis about no parameter change occurring at the point $\tau$, we take an interval $\mathcal{I}'$ of roughly the same length as $\mathcal{I}''$ and use as a test statistic the difference between the sum of log likelihoods corresponding to the intervals $\mathcal{I}'$ and $\mathcal{I}''$, and the log likelihood corresponding to the interval $\mathcal{I} = \mathcal{I}' \cup \mathcal{I}''$ containing no change points:

$$T_{\mathcal{I},\tau} = \max_{\theta',\theta''}(L_{\mathcal{I}''}(\theta'') + L_{\mathcal{I}'}(\theta')) - \max_{\theta} L_{\mathcal{I}}(\theta) = L_{\mathcal{I}'}(\tilde{\theta}_{\mathcal{I}'}) + L_{\mathcal{I}''}(\tilde{\theta}_{\mathcal{I}''}) - L_{\mathcal{I}}(\tilde{\theta}_{\mathcal{I}}), \tag{7}$$

where $L(\cdot)$ denotes the log likelihood function. Under the null hypothesis the test statistic should not exceed a certain *critical value* $\mathfrak{z}$. Failure to reject the null hypothesis implies that $\tau$ is not a change point, hence the interval $\mathcal{I}''$ can be extended by including $\tau$.

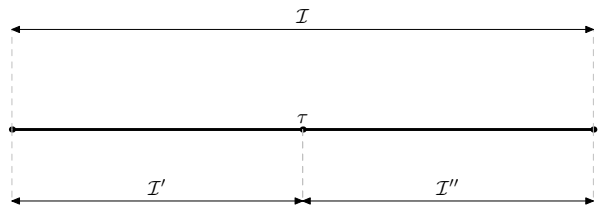

FIG 2. *Intervals involved in the change point detection procedure.*

To test the interval $I_k$ provided $I_{k-1}$ contains no change points, we let $\mathcal{I} = I_{k+1}$ and accept $I_k$ if the above null hypothesis cannot be rejected for every point $\tau \in I_k \setminus I_{k-1}$ (Fig. 2). Therefore, the predicate $\mathcal{A}$ from (6) assumes the form

$$\mathcal{A}(k) = \mathbf{1}\{T_{I_{k+1},\tau} \leq \mathfrak{z}_k \quad \text{for each } \tau \in I_k \setminus I_{k-1}\}.$$

Since an enclosing interval is necessary for testing, the largest interval that can be tested by LCP is $I_{K-1}$.

## 2.2. Local model selection

The local model selection procedure goes back to Lepski (1990), see also Spokoiny and Vial (2008). Given a sequence of accepted weak estimates $\tilde{\theta}_1 \ldots \tilde{\theta}_{k-1}$, the procedure tests a candidate estimate $\tilde{\theta}_k$. The candidate estimate is accepted if it belongs to the condifence interval $\mathcal{E}$ of each of the previous weak estimates (cf. Figure 3). Formally, the predicate $\mathcal{A}$ for LMS has the form

$$\mathcal{A}(k) = \mathbf{1}\{\tilde{\theta}_k \in \mathcal{E}_l \text{ for all } l < k\}.$$

The procedure for constructing the confidence intervals of the weak estimates is provided by the corollary of the Theorem 2 that claims $\mathcal{E}(\mathfrak{z}_\alpha) = \left\{\theta : L(\tilde{\theta},\theta) \leq \mathfrak{z}_\alpha\right\}$ to be the $\alpha$-confidence set for the parameter $\theta$.

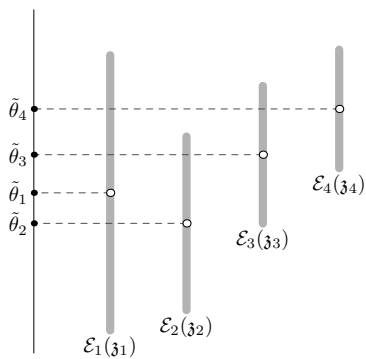

FIG 3. *Principle of the local model selection.* $\hat{\theta} = \tilde{\theta}_3$.

## 2.3. Stagewise aggregation

The SA procedure introduced in Belomestny and Spokoiny (2007) differs from the two methods described above in that it does not choose the adaptive estimate $\hat{\theta}$ from the set of weak estimates $\tilde{\theta}_1 \ldots \tilde{\theta}_K$. Instead, based on the weak estimates, it sequentially constructs *aggregated estimates* $\hat{\theta}_1 \ldots \hat{\theta}_K$ possessing the property that any aggregated estimate $\hat{\theta}_k$ has smaller variance than the corresponding weak estimate $\tilde{\theta}_k$, while keeping "close" to it in terms of the statistical difference, the latter being measured through the fitted likelihood $L_k(\tilde{\theta}_k, \hat{\theta}_{k-1}) = L_k(\tilde{\theta}_k) - L_k(\hat{\theta}_{k-1})$. The adaptive estimate is finally taken equal to the last aggregated estimate: $\hat{\theta} = \hat{\theta}_K$, unless the *early stopping* (see below) occurs.

Formally, the first aggregated estimate is equal to the first weak estimate and every next aggregated estimate is a convex combination of the previous aggregated estimate and the current weak estimate:

$$\hat{\theta}_k = \begin{cases} \tilde{\theta}_1, & k = 1 \\ \gamma_k \tilde{\theta}_k + (1 - \gamma_k) \hat{\theta}_{k-1}, & k = 2, \ldots, K \end{cases} \tag{8}$$

Here $\gamma_k$ is the *mixing coefficient* that reflects the statistical difference between the previous aggregated estimate $\hat{\theta}_{k-1}$ and the current weak estimate $\tilde{\theta}_k$, and is obtained by applying an *aggregation kernel* $K_{\mathrm{ag}}$ to the fitted likelihood $L_k(\tilde{\theta}_k, \hat{\theta}_{k-1})$ scaled by the *critical value* $\mathfrak{z}_k$:

$$\gamma_k = K_{\mathrm{ag}}\left(\frac{L_k(\tilde{\theta}_k, \hat{\theta}_{k-1})}{\mathfrak{z}_k}\right).$$

The aggregation kernel acts as a link between the likelihood ratio and the mixing coefficient. The kernel is selected to ensure that a smaller statistical difference between $\tilde{\theta}_k$ and $\hat{\theta}_{k-1}$ leads to the mixing coefficient close to 1 and thus to the aggregated estimate $\hat{\theta}_k$ close to $\tilde{\theta}_k$, while a large difference should provide the mixing coefficient close to zero and thus keep $\hat{\theta}_k$ close to $\hat{\theta}_{k-1}$.

Whenever the difference is very large, the mixing coefficient is zero, and the procedure stops prematurely by setting $\hat{\theta} = \hat{\theta}_{k-1}$. We call this situation *early stopping*.

To satisfy the mentioned requirements, the kernel must be supported on the closed interval $[0, 1]$ and monotonously decrease from 1 on the left edge to 0 on the right edge and have a plateau of size $b$ at the left. Thus, the aggregation kernel assumes the form:

$$K_{\mathrm{ag}}(u) = \begin{cases} 1, & 0 \le u < b \\ 1 - \bar{K}_{\mathrm{ag}}(u), & b \le u \le 1 \end{cases}.$$

Examples of $\bar{K}_{\mathrm{ag}}(u)$ include $\frac{u-b}{1-b}$ (triangular kernel), $\left(\frac{u-b}{1-b}\right)^2$ (Epanechnikov kernel) etc.

## 3. Critical values and other parameters

All procedures described above depend on a set of parameters $\mathfrak{z}_1 \ldots \mathfrak{z}_K$ that we call critical values. The critical values reflect the problem design (interval length, model, method etc.). The principle behind their selection is to ensure that the probability of error of the first kind under the null hypothesis (in the parametric case) does not exceed a prescribed level, while errors on initial estimation steps (for small $k$) are penalized stronger. To formalize this idea, we introduce the notion of the adaptive estimate obtained on the $k$-th step $\hat{\theta}_k$ for model selection-type procedures (for the case of SSA see (8)):

$$\hat{\theta}_k = \begin{cases} \tilde{\theta}_k, & \mathcal{A}(l) = 1 \quad \text{for all } l \le k, \\ \hat{\theta}_{k-1} & \text{otherwise.} \end{cases}$$

The selection of the critical values is based on the following *propagation condition*, that postulates the performance in the parametric case:

**Condition 2 (Propagation).** *Let the observations be i.i.d. from $\mathbf{P}_{\theta^*}$ with $\theta^* \in \Theta$ and $\alpha$ be a constant, $0 < \alpha \le 1$. Then the risk of the adaptive estimate in the parametric case $\hat{\theta}_k$ attained on the $k$-th step must be bounded:*

$$\frac{\mathbf{E}_{\theta^*} |L_k(\tilde{\theta}_k, \hat{\theta}_k)|^r}{\mathcal{R}_{r,\theta^*,k}} \le \alpha \frac{k}{K} \quad \text{for } k = 1, \ldots, K, \tag{9}$$

where $\mathcal{R}_{r,\theta^*,k}$ is the parametric risk defined analogously to (5) with $I = I_k$.

In the parametric case the adaptive estimate on the $k$-th step should be equal to the oracle estimate: $\hat{\theta}_k = \tilde{\theta}_k$. Otherwise the *false alarm* has occurred, incurring the loss $L_k(\tilde{\theta}_k, \hat{\theta}_k)$. The propagation condition stipulates that the risk (mean loss) associated with the $k$-th adaptive estimate does not exceed a certain fraction of the parametric risk. The constant $\alpha$ plays the role of the level of the procedure.

However, the propagation condition is not explicit. We use the following sequential method for the computation of critical values using Monte-Carlo simulations. Denote as $\hat{\theta}_l(\mathfrak{z}_k)$ for $l \geq k$ the adaptive estimate obtained after the $l$-th step of the procedure run with the critical values $\mathfrak{z}_1, \ldots, \mathfrak{z}_{k-1}$ known and $\mathfrak{z}_{k+1}, \ldots, \mathfrak{z}_K$ set to infinity:

$$\hat{\theta}_l(\mathfrak{z}_k) = \hat{\theta}_l(\mathfrak{z}_1, \ldots, \mathfrak{z}_k, \mathfrak{z}_{k+1} = \infty, \ldots, \mathfrak{z}_K = \infty), \quad l > k.$$

The first critical value can be selected to satisfy the conditions

$$\frac{\mathbf{E}_{\theta^*} \left| L_l(\tilde{\theta}_l, \hat{\theta}_l(\mathfrak{z}_1)) \right|^r}{\mathcal{R}_{r,\theta^*,l}} \leq \frac{\alpha}{K}, \quad l = 2, \ldots, K.$$

Such a value exists, since for $\mathfrak{z}_1$ taken sufficiently large the weak and adaptive estimates coincide for any $l$ and all Monte-Carlo paths from $\mathbf{P}_{\theta^*}$, leading to the zero risk. With the first $k-1$ critical values fixed, the $k$-th critical value is selected using the condition

$$\frac{\mathbf{E}_{\theta^*} \left| L_l(\tilde{\theta}_l, \hat{\theta}_l(\mathfrak{z}_k)) \right|^r}{\mathcal{R}_{r,\theta^*,l}} \leq \alpha \frac{k}{K}, \quad l = k+1, \ldots, K.$$

Obtaining critical values is computationally involved, however for a given parametric setup it is sufficient to compute them once under the parametric model and reuse later. Critical values for various models are given in the Section 5.1, where the influence of the parameters $\alpha$ (test level) and $r$ is discussed as well.

## 4. Theoretical properties of adaptive estimates

In this section the estimates obtained using the methods described in the previous sections are shown to possess some remarkable properties. Namely,

(a) in the local parametric case (under the small modelling bias condition) the difference between the oracle estimate and the adaptive estimate on each step is only within a small factor of the corresponding difference in the global parametric case, and the factor is due to the "payment for adaptation" (propagation result);
(b) the quality of the adaptive estimates does not deteriorate upon violation of the small modelling bias condition (stability result);
(c) the risk between the final adaptive estimate and the oracle estimate is bounded (oracle result).

Some of these results have been obtained previously for particular methods, e.g. in Katkovnik and Spokoiny (2008) for LMS, Mercurio and Spokoiny (2004) for LCP and Belomestny and Spokoiny (2007) for SSA. Here they are reformulated for the general case of any of the three procedures.

We introduce two regularity conditions needed for the formulation of the claimed properties.

**Condition 3 (Exponential growth of the intervals).** *For some constants $\mathfrak{u}_0, \mathfrak{u}$ s. t. $\mathfrak{u}_0 \leq \mathfrak{u} < 1$ the interval lengths $N_k$ satisfy for every $2 \leq k \leq K$*

$$\mathfrak{u}_0 \leq \frac{N_{k-1}}{N_k} \leq \mathfrak{u}.$$

**Condition 4 (Compactness of the parameter set).** *There exists a value $\mathfrak{a} \geq 1$, s. t. for any two parameter values $\theta_1, \theta_2 \in \Theta$*

$$\frac{1}{\mathfrak{a}^2} \leq \frac{d''(\theta_1)}{d''(\theta_2)} \leq \mathfrak{a}^2,$$

*where d is the function from the definition of the exponential family, see Section 1.2.*

**Definition (Oracle index).** *An* oracle index *is the number $k^\circ$ of the largest interval satisfying Condition 1 with $I^\circ = I_{k^\circ}$ for some $\theta^\circ \in \Theta$ and $\Delta^\circ \geq 0$.*

### 4.1. Propagation property

The following result, implied by the Theorem 4 and Condition 2, states that under the small modelling bias condition the performance of the adaptive procedures is essentially the same as in the parametric case.

**Theorem 5 (Propagation property).** *Assume the regularity Conditions 3 and 4. Then under Condition 1 for any $r \geq 0$*

$$\mathbf{E}\log\left(1 + \frac{\left|L_{k^\circ}(\tilde{\theta}_{k^\circ}, \hat{\theta}_{k^\circ})\right|^r}{\mathcal{R}_{r,\theta^\circ,k^\circ}}\right) \leq \alpha + \Delta^\circ$$

*and*

$$\mathbf{E}\log\left(1 + \frac{\left|L_{k^\circ}(\tilde{\theta}_{k^\circ}, \theta^\circ)\right|^r}{\mathcal{R}_{r,\theta^\circ,k^\circ}}\right) \leq 1 + \Delta^\circ.$$

### 4.2. Stability after propagation

The next result, following from the definition of the LMS procedure, claims that the quality of the adaptive estimate does not deteriorate with the growing interval even when the small modelling bias condition does not hold anymore.

**Theorem 6 (Stability after propagation for LMS).** *For $\hat{k} \geq k^\circ$, where $k^\circ$ is the oracle index, the following holds:*

$$L_{k^\circ}(\tilde{\theta}_{k^\circ}, \hat{\theta}) = L_{k^\circ}(\tilde{\theta}_{k^\circ}, \tilde{\theta}_{\hat{k}}) \leq \mathfrak{z}_{k^\circ}. \tag{10}$$

Analogous result for the LCP method is established in the Theorem 7.

**Theorem 7 (Stability for LCP).** *It holds for every $k < K$*

$$N_k \mathcal{K}(\hat{\theta}_k, \hat{\theta}_{k+1}) \leq \mathfrak{z}_k. \tag{11}$$

*Moreover, if Conditions 3 and 4 hold, then for every $k < k' \leq K$*

$$N_k \mathcal{K}(\hat{\theta}_k, \hat{\theta}_{k'}) \leq \left(\frac{\mathfrak{a}}{1-\sqrt{\mathfrak{u}}}\right)^2 \bar{\mathfrak{z}}_k, \tag{12}$$

*where $\bar{\mathfrak{z}}_k = \max_{l \geq k} \mathfrak{z}_l$.*

**Proof.** If $I_{k+1}$ is rejected, then $\hat{\theta}_{k+1} = \hat{\theta}_k$ and the assertion (11) trivially follows. Whenever $I_{k+1}$ is accepted, it holds that $\hat{\theta}_k = \tilde{\theta}_k$ and $\hat{\theta}_{k+1} = \tilde{\theta}_{k+1}$. The acceptance of $I_k$ implies by definition of the procedure that $T_{I_k} \leq \mathfrak{z}_k$ and, in particular, $T_{I_k, \tau} \leq \mathfrak{z}_k$ with $\tau = t^\diamond - N_k$ being the left edge of $I_k$. Due to the definition (7) of the test statistic and taking into account the representation (4) of the fitted likelihood for exponential families it follows that

$$L_k(\tilde{\theta}_k, \tilde{\theta}_{k+1}) = N_k \mathcal{K}(\tilde{\theta}_k, \tilde{\theta}_{k+1}) \leq \mathfrak{z}_k.$$

That proves the assertion (11).

Decomposing the Kullback – Leibler divergence according to Lemma 3 (which requires Condition 4) and applying the just proven assertion (11) yield

$$\mathcal{K}^{1/2}(\hat{\theta}_k, \hat{\theta}_{k'}) \leq \mathfrak{a} \sum_{j=k}^{k'-1} \mathcal{K}^{1/2}(\hat{\theta}_j, \hat{\theta}_{j+1}) \leq \mathfrak{a} \sum_{j=k}^{k'-1} \sqrt{\frac{\mathfrak{z}_j}{N_j}}.$$

Use of Condition 3 leads to the bound

$$\mathcal{K}^{1/2}(\hat{\theta}_k, \hat{\theta}_{k'}) \leq \mathfrak{a}\sqrt{\frac{\bar{\mathfrak{z}}_k}{N_k}} \sum_{j=k}^{k'-1} \left(\frac{N_k}{N_{k+1}} \frac{N_{k+1}}{N_{k+2}} \cdots \frac{N_{j-1}}{N_j}\right)^{1/2} \leq \mathfrak{a}\sqrt{\frac{\bar{\mathfrak{z}}_k}{N_k}} \sum_{j=k}^{k'-1} \mathfrak{u}^{(j-k)/2}.$$

The last expression is a partial sum of the geometric series with common ratio $\sqrt{\mathfrak{u}} < 1$. Therefore, it can be bounded from above by the total sum:

$$\mathfrak{a}\sqrt{\frac{\bar{\mathfrak{z}}_k}{N_k}} \sum_{j=k}^{k'-1} \mathfrak{u}^{(j-k)/2} \leq \frac{\mathfrak{a}}{1-\sqrt{\mathfrak{u}}} \sqrt{\frac{\bar{\mathfrak{z}}_k}{N_k}}.$$

Consequently,

$$\left(N_k \mathcal{K}(\hat{\theta}_k, \hat{\theta}_{k'})\right)^{1/2} \leq \frac{\mathfrak{a}}{1-\sqrt{\mathfrak{u}}} \sqrt{\bar{\mathfrak{z}}_k}.$$

Squaring both sides leads to the assertion (12). □

The stability result for the stagewise aggregation method is formulated in the Theorem 8.

**Theorem 8 (Stability for the SSA).** *It holds for every $k \leq K$*

$$N_k \mathcal{K}(\hat{\theta}_k, \hat{\theta}_{k-1}) \leq \mathfrak{z}_k. \tag{13}$$

*Moreover, if Condition 3 holds, then for every $k < k' \leq K$*

$$N_k \mathcal{K}(\hat{\theta}_k, \hat{\theta}_{k'}) \leq \mathfrak{a}^2 c_{\mathfrak{u}}^2 \bar{\mathfrak{z}}_k, \tag{14}$$

*where $c_{\mathfrak{u}} = 1/(1 - \sqrt{\mathfrak{u}})$ and $\bar{\mathfrak{z}}_k$ has been defined in Theorem 7.*

**Proof.** Convexity of the Kullback – Leibler divergence $\mathcal{K}(u, v)$ w.r.t. the first argument implies

$$\mathcal{K}(\hat{\theta}_k, \hat{\theta}_{k-1}) = \mathcal{K}(\gamma_k \tilde{\theta}_k + (1 - \gamma_k)\hat{\theta}_{k-1}, \hat{\theta}_{k-1})$$
$$\leq \gamma_k \mathcal{K}(\tilde{\theta}_k, \hat{\theta}_{k-1}) + (1 - \gamma_k)\mathcal{K}(\hat{\theta}_{k-1}, \hat{\theta}_{k-1}) = \gamma_k \mathcal{K}(\tilde{\theta}_k, \hat{\theta}_{k-1}).$$

If $\mathcal{K}(\hat{\theta}_k, \hat{\theta}_{k-1}) \geq \mathfrak{z}_k/N_k$, then $\gamma_k = 0$ and (13) follows. Assertion (14) is proven exactly as the analogous assertion (12). □

### 4.3. Oracle result

Finally, we show that the quality of adaptive estimates is comparable with that of the oracle estimate. In the following theorems we assume the small modelling bias condition (Condition 1) as well as regularity conditions (Conditions 3 and 4) to hold.

We begin with the oracle property of the LMS estimates.

**Theorem 9.** *For any $r \geq 0$ the LMS adaptive estimate $\hat{\theta}$ satisfies*

$$\mathbf{E} \log \left(1 + \frac{\left|L_{k^\circ}(\tilde{\theta}_{k^\circ}, \hat{\theta})\right|^r}{\mathfrak{z}_{k^\circ}^r}\right) \leq \Delta^\circ + \alpha \mathfrak{r}_r / \mathfrak{z}_{k^\circ}^r + 1.$$

**Proof.** Definition of the LMS adaptive estimate $\hat{\theta} = \tilde{\theta}_{\hat{k}}$ and the stability property (10) imply

$$\left|L_{k^\circ}(\tilde{\theta}_{k^\circ}, \hat{\theta})\right|^r = \left|L_{k^\circ}(\tilde{\theta}_{k^\circ}, \hat{\theta}_{k^\circ})\right|^r + \left|L_{k^\circ}(\tilde{\theta}_{k^\circ}, \tilde{\theta}_{k^\circ})\right|^r \mathbf{1}\{\hat{k} > k^\circ\} \leq \left|L_{k^\circ}(\tilde{\theta}_{k^\circ}, \hat{\theta}_{k^\circ})\right|^r + \mathfrak{z}_{k^\circ}^r. \tag{15}$$

By Condition 2 and Theorem 3

$$\mathbf{E}_{\theta^\circ} \left|L_{k^\circ}(\tilde{\theta}_{k^\circ}, \hat{\theta}_{k^\circ})\right|^r \leq \alpha \mathcal{R}_{r,\theta^\circ} \leq \alpha \mathfrak{r}_r.$$

By Lemma 2

$$\mathbf{E} \log \left(1 + \frac{\left|L_{k^\circ}(\tilde{\theta}_{k^\circ}, \hat{\theta})\right|^r}{\mathfrak{z}_{k^\circ}^r}\right) \leq \Delta^\circ + \frac{\mathbf{E}_{\theta^\circ}\left|L_{k^\circ}(\tilde{\theta}_{k^\circ}, \hat{\theta})\right|^r}{\mathfrak{z}_{k^\circ}^r} \leq \Delta^\circ + \frac{\alpha \mathfrak{r}_r}{\mathfrak{z}_{k^\circ}^r} + 1,$$

and the assertion follows. □

Since the decomposition (15) is only possible for the LMS estimates, the oracle result for the LCP and SSA estimates is somewhat weaker. It is formulated in the following theorem.

**Theorem 10.** *The estimate $\hat{\theta}$ is close to the oracle estimate $\tilde{\theta}_{k^\circ}$ in the sense that*

$$\mathbf{E}\log\left(1+\frac{\left|L_{k^\circ}(\tilde{\theta}_{k^\circ},\hat{\theta}_{k^\circ})\right|^r}{\mathfrak{r}_r}\right) \leq \alpha + \Delta^\circ,$$

$$N_{k^\circ}\mathcal{K}(\hat{\theta}_{k^\circ},\hat{\theta}) \leq \mathfrak{a}^2 c_\mathfrak{u}^2 \bar{\mathfrak{z}}_{k^\circ},$$

*where $\mathfrak{r}_r$ is the parametric risk bound from the Theorem 3.*

We prove the corollary for the case of $r = 1/2$.

**Theorem 11.** *Assume regularity Conditions 3 and 4. Then under Condition 1 the following holds:*

$$\mathbf{E}\log\left(1+\frac{\left|N_{k^\circ}\mathcal{K}(\hat{\theta},\theta^\circ)\right|^{1/2}}{\mathfrak{a}\mathfrak{r}_{1/2}}\right) \leq \log\left(1+\frac{c_\mathfrak{u}\sqrt{\bar{\mathfrak{z}}_{k^\circ}}}{\mathfrak{r}_{1/2}}\right) + \Delta^\circ + \alpha + 1, \qquad (16)$$

*where $c_\mathfrak{u}$ and $\bar{\mathfrak{z}}_k$ has been defined in Theorem 8.*

**Proof.** By Lemma 3 and stability property (12)

$$\frac{\left|N_{k^\circ}\mathcal{K}(\hat{\theta},\theta^\circ)\right|^{1/2}}{\mathfrak{a}} \leq \left|N_{k^\circ}\mathcal{K}(\hat{\theta}_{k^\circ},\hat{\theta})\right|^{1/2} + \left|N_{k^\circ}\mathcal{K}(\tilde{\theta}_{k^\circ},\hat{\theta}_{k^\circ})\right|^{1/2} + \left|N_{k^\circ}\mathcal{K}(\tilde{\theta}_{k^\circ},\theta^\circ)\right|^{1/2}$$

$$\leq c_\mathfrak{u}\sqrt{\bar{\mathfrak{z}}_{k^\circ}} + \left|L_{k^\circ}(\tilde{\theta}_{k^\circ},\hat{\theta}_{k^\circ})\right|^{1/2} + \left|L_{k^\circ}(\tilde{\theta}_{k^\circ},\theta^\circ)\right|^{1/2}.$$
(17)

Further, it is true that for all $a, b \geq 0$

$$\log(1+a+b) \leq \log(1+a) + \log(1+b). \qquad (18)$$

Substituting (17) into (16) and applying the inequality (18) yields:

$$\mathbf{E}\log\left(1+\frac{\left|N_{k^\circ}\mathcal{K}(\hat{\theta},\theta^\circ)\right|^{1/2}}{\mathfrak{a}\mathfrak{r}_{1/2}}\right) \leq \log\left(1+\frac{c_\mathfrak{u}\sqrt{\bar{\mathfrak{z}}_{k^\circ}}}{\mathfrak{r}_{1/2}}\right) + \mathbf{E}\log(1+\zeta)$$

with $\zeta = \mathfrak{r}_{1/2}^{-1}\left[\left|L_{k^\circ}(\tilde{\theta}_{k^\circ},\hat{\theta}_{k^\circ})\right|^{1/2} + \left|L_{k^\circ}(\tilde{\theta}_{k^\circ},\theta^\circ)\right|^{1/2}\right]$. Applying Lemma 2 and estimating $\mathbf{E}_{\theta^\circ}\zeta$ in view of the propagation condition (9) and Theorem 3 complete the proof. □

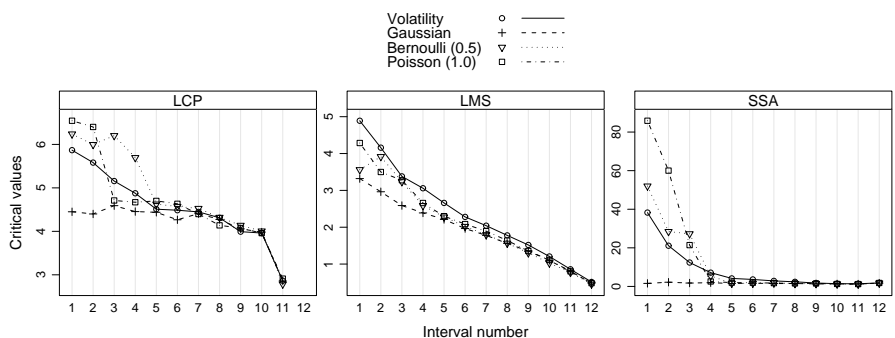

FIG 4. *Critical values for the standard setup (see text). True parameter value, if relevant, is indicated in parentheses after the model designation.*

## 5. Simulations and applications

### 5.1. Influence of the parameters on the critical values

From the propagation condition (Condition 2) is follows that the critical values depend on global parameters $\alpha$ (level of the procedure) and $r$ as well as on the underlying family of distributions and on the sequence on intervals generating weak estimates. In this section the influence of the mentioned parameters on the critical values is investigated.

Figure 4 shows the critical values for the standard setup: $\alpha = 0.7$, $r = 0.5$, interval length from 5 to 284 points growing exponentially with the factor 1.4. Other figures in the section show the critical values for the setups other than the standard one. In every case only the indicated parameter is being changed, the other remaining at the standard settings. For the Gaussian and volatility models the critical values do not depend on the unknown true parameter $\theta^*$ as stated by the following lemma.

**Lemma 1 (Pivotality of the critical values).** *Let the observations $Y_t$ be i.i.d. from $\mathbf{P}_{\theta^*}$ with $\theta^* \in \Theta$, where $\mathbf{P}$ describes either Gaussian or volatility model. Then the critical values obtained using the propagation condition (Condition 2) do not depend on $\theta^*$.*

**Proof.** It suffices to show the pivotality of the fitted likelihood $L_k(\tilde{\theta}_k, \hat{\theta}_k) = N_k \mathcal{K}(\tilde{\theta}_k, \hat{\theta}_k)$ entering Condition 2. In the rest of the proof we drop the index $k$. Observations originating from the Gaussian model are of the form

$$Y_t = \theta^* + \sigma \varepsilon_t, \quad \varepsilon_t \sim \mathbb{N}(0, 1),$$

where $\sigma$ is a known constant. Hence, weak estimates assume the representation

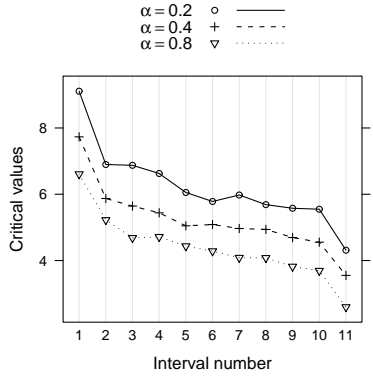

Fig 5. *Critical values corresponding to various test levels (model: Volatility, method: LCP, standard intervals).*

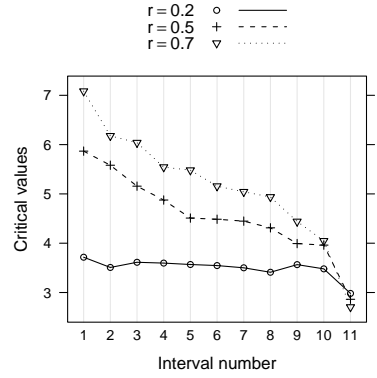

Fig 6. *Critical values corresponding to various r levels (model: Volatility, method: LCP, standard intervals).*

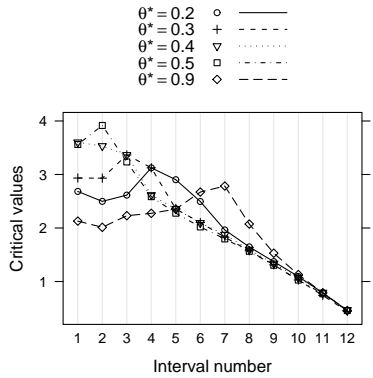

Fig 7. *Critical values for various $\theta^*$ (model: Bernoulli, method: LMS, standard parameters and intervals).*

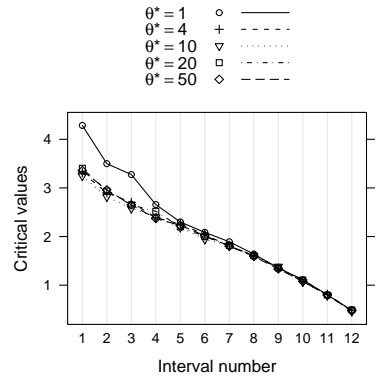

Fig 8. *Critical values for various $\theta^*$ (model: Poisson, method: LMS, standard parameters and intervals).*

as a sum of $\theta^*$ and a further item independent of $\theta^*$:

$$\tilde{\theta} = \frac{1}{N}\sum_{t=1}^{N} Y_t = \theta^* + \frac{\sigma}{N}\sum_{t=1}^{N}\varepsilon_t,$$

This is also true for the adaptive estimate $\hat{\theta}$ by construction of the procedures LCP and LMS and can be shown to hold for SSA by induction. Since the Kullback – Leibler divergence $\mathcal{K}(\theta_1, \theta_2)$ for the Gaussian model depends only on the difference of the arguments (cf. Table 1), the true parameter value $\theta^*$ cancels.

Observations following the volatility model assume the form

$$Y_t = \theta^*\varepsilon_t^2, \quad \varepsilon_t \sim \mathbb{N}(0,1),$$

and the weak estimates can be represented as a product of the true parameter value and a further term:

$$\tilde{\theta} = \frac{\theta^*}{N}\sum_{t=1}^{N}\varepsilon_t^2,$$

with the same holding for adaptive estimates. As the Kullback – Leibler divergence $\mathcal{K}(\theta_1, \theta_2)$ for the volatility model is a function of the ratio of the arguments (cf. Table 1), $\theta^*$ cancels. This concludes the proof. □

On the Figures 7 and 8 several sets of the critical values corresponding to various values of $\theta^*$ are presented. Dependence of the critical values on $\theta^*$ is not very strong and it is sufficient in practice to compute one set of critical values for a certain $\theta^*$, e.g. $\theta^* = 0.5$ in the case of Bernoulli model.

### 5.2. Delay of detection of a jump

A jump in the the parameter of the underlying distribution cannot be detected immediately, but after a delay at least as large as the smallest interval. The subject of this section is to investigate this delay. Consider a sequence of the true parameter values $\theta^*$ exhibiting a jump of height $h$ from a value $\underline{\theta}$ to a value $\overline{\theta} = \underline{\theta} + h$ at time $t^\diamond$ (Figure 9). We regard the jump as detected once the parameter estimate has reached the $\delta$-fraction of the difference between $\overline{\theta}$ and $\underline{\theta}$ for some $0 < \delta \leq 1$:

$$\mathfrak{t} = \min\{t: \quad \hat{\theta}_t \geq \underline{\theta} + \delta(\overline{\theta} - \underline{\theta})\}.$$

We present the mean delay $\mathfrak{t} - t^\diamond$ of 1000 replications (Figure 10) as a function of $\overline{\theta}$ for various distributions and $\delta = 0.7$ . As one would await, smaller jumps require more time to be detected.

### 5.3. Switching regime model

In the present section adaptive methods are applied to data following a switching regime model. The latter is described by a sequence of Markov moments $\{\nu_i\} = \{40, 99, 85, 99, 51, 38\}$ w.r.t. the filtration $(\mathcal{F}_t)$, and a sequence of states $\{\theta_{n_i}\}$ with $\{n_i\} = \{2, 4, 3, 5, 2, 1\}$. The states for various models are listed in the Table 2. Figure 11 illustrates estimation by different methods based on 10000 realizations of data originating from four distributions. One can recognize that the methods in use demonstrate very reasonable performance. The adaptive estimate is most of the time within the oracle confidence bounds. Note that oracle confidence bounds reflect the jumps of the true value with some delay, leading to skewness of the "wavefronts". This is due to the fact that by construction of the estimation procedures the smallest interval is always accepted, delaying the oracle.

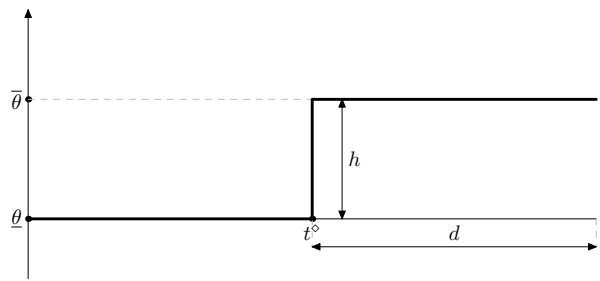

Fig 9. *Design of the jump experiment.*

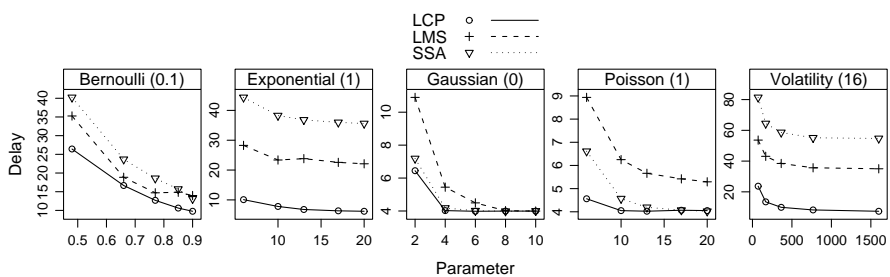

Fig 10. *Mean delay of detecting a jump for various distributions (the number in parentheses is $\underline{\theta}$ ), $\delta = 0.7$.*

The estimation is characterized numerically in terms of the mean squared error, Kullback–Leibler divergence and mean absolute error in the Table 3.

Performance of the LCP and LMS is very close, while the SA works somewhat worse on data with abrupt jumps.

### 5.4. Application to volatility estimation

In this section the performance of the procedures is demonstrated by means of applying them to the problem of volatility estimation of the daily stock price of Allianz in the period from 1974-01-02 to 1996-12-30.

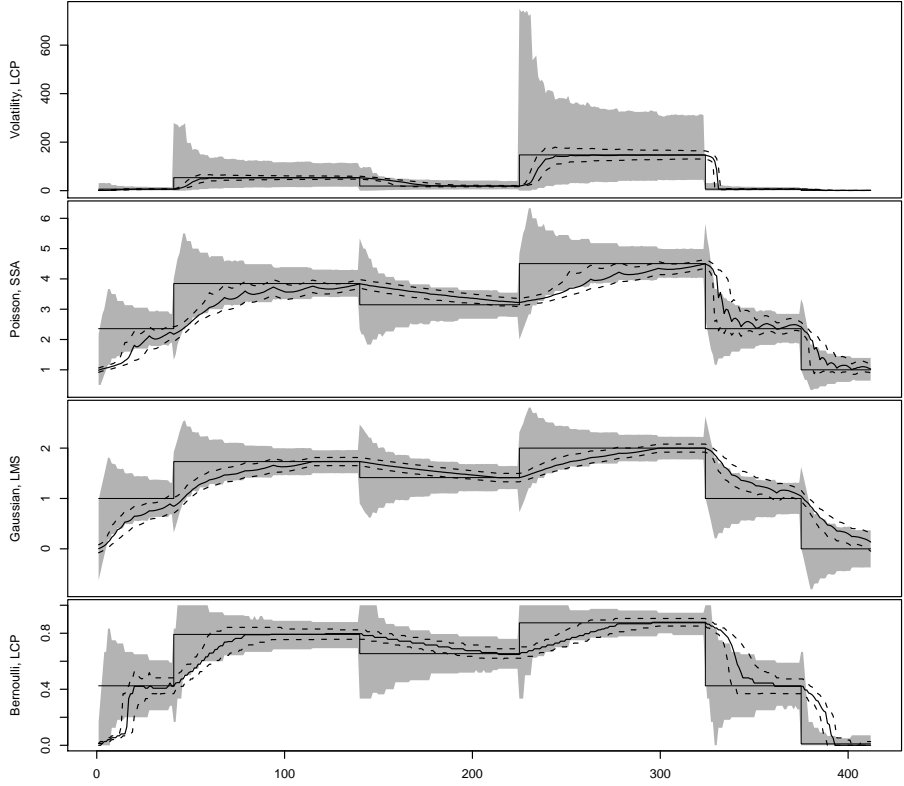

FIG 11. *Estimation for selected combinations of method and model. Thin solid line: true value; shaded area: 95% confidence region for the oracle estimate; thick solid line: median estimate; dashed lines: 1st and 3rd quartiles.*

| Model | 1 | 2 | 3 | 4 | 5 |
|---|---|---|---|---|---|
| Gaussian | 0.00 | 1.00 | 1.41 | 1.73 | 2.00 |
| Poisson | 1.00 | 2.36 | 3.15 | 3.85 | 4.51 |
| Bernoulli | 0.01 | 0.42 | 0.65 | 0.79 | 0.87 |
| Volatility | 1.00 | 6.31 | 19.06 | 53.59 | 147.41 |

TABLE 2
*States of the switching regime model for various distributions.*

| | Method | MSqE | KLD | MAE |
|---|---|---|---|---|
| **Volatility** | | | | |
| | LCP | 567.66 | 0.07 | 9.08 |
| | LMS | 674.27 | 0.08 | 12.02 |
| | SSA | 1101.21 | 0.14 | 16.48 |
| **Gaussian** | | | | |
| | LCP | 0.08 | 0.04 | 0.17 |
| | LMS | 0.07 | 0.03 | 0.18 |
| | SSA | 0.09 | 0.04 | 0.20 |
| **Poisson** | | | | |
| | LCP | 0.23 | 0.04 | 0.29 |
| | LMS | 0.20 | 0.04 | 0.31 |
| | SSA | 0.27 | 0.05 | 0.35 |
| **Bernoulli** | | | | |
| | LCP | 0.02 | 0.11 | 0.08 |
| | LMS | 0.01 | 0.05 | 0.08 |
| | SSA | 0.04 | 0.15 | 0.14 |

TABLE 3
*Estimation quality in terms of the mean squared error, Kullback–Leibler divergence and mean absolute error. Based on 10000 replications.*

Let $S_t$ be an observed stock price process, and $Y_t = \left(\log \frac{S_t}{S_{t-1}}\right)^2$ the squared log returns (Figure 12). The latter are described by the *conditional heteroskedastic* model

$$Y_t = \theta_t \varepsilon_t^2, \quad \varepsilon_t \sim \mathbb{N}(0,1),$$

which is a particular case of the general model of Section 1. The problem is to forecast $\theta_{t+h}$ from $Y_1, \ldots, Y_t$ for the time horizon $h$. The estimation is carried out using the procedures introduced above. The criterion to describe and compare the performance of the procedures is the median $h$-step-ahead forecasting error defined as

$$\text{FE}_h = \text{median}[\mathcal{K}(\bar{Y}_{t,h}, \hat{\theta}_t)],$$

where $\bar{Y}_{t,h} = \sum_{h' \leq h} Y_{t+h'}/h$ and the horizon $h$ is taken equal to 1, 3 and 10 steps. The results are shown in the Figure 13. All procedures demonstrate comparable and very reasonable performance. The stagewise aggregation procedure performs somewhat better than the other methods, in contrast to the switching regime model of Section 5.3.

## 5.5. Application to waiting time estimation

With the advent of fast and cheap computers equipped with large storage high frequency financial data became available, opening new perspectives for financial time series analysis such as the study of bid-ask spread, transaction intensity, waiting times etc. In this section we apply adaptive methods to the estimation of waiting times, i.e. time intervals between transactions. Waiting times are assumed to follow the exponential distribution with the density

$$p(x, \theta) = \frac{1}{\theta} e^{-x/\theta}.$$

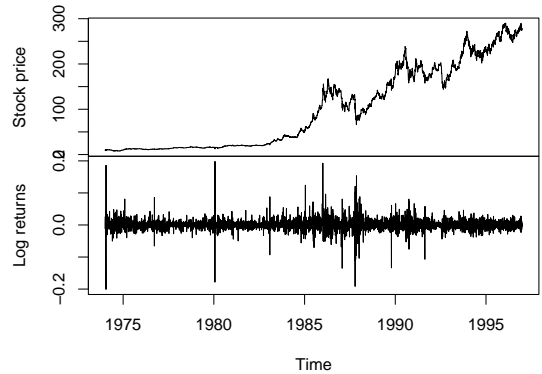

Fig 12. *Price and log returns of the Allianz stock.*

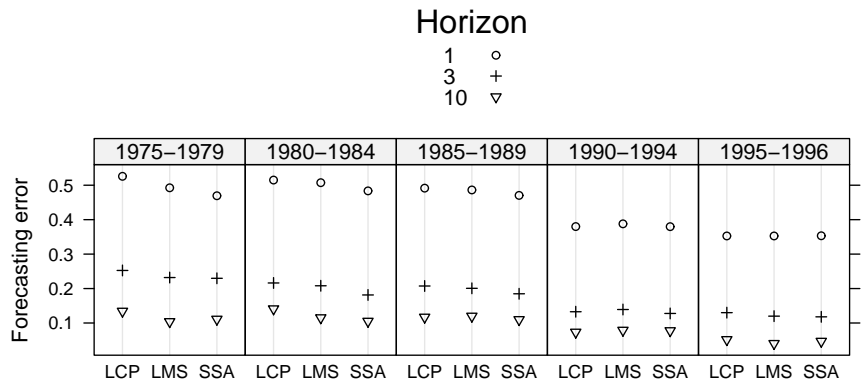

Fig 13. *Forecasting error of the Allianz stock price volatility estimation.*

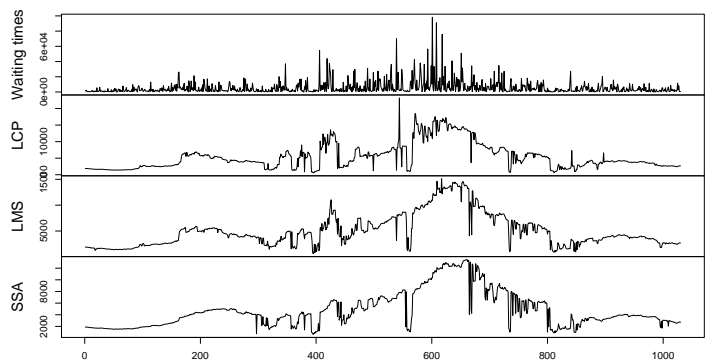

FIG 14. *Upper panel: waiting times between transactions (BHP Billiton stock on the Australian Stock Exchange on July 5, 2002); lower panels: corresponding adaptive estimates.*

For an example we use the data on the BHP Billiton stock traded on the Australian stock exchange (ASX) from July 1, 2002 to August 30, 2002. Figure 14 displays waiting times within one day as well as the estimates of $\theta$. Time intervals are measured in ticks with the tick size of 0.00390625 seconds. As in the section 5.4, the median forecasting error is used as the performance criterion with the horizon $h = 10$. The forecasting error of the three methods are presented in the Figure 15. The LMS method appears to perform best in this case.

**Conclusion**

In the paper a unified approach to the adaptive estimation of the univariate time series parameters have been presented and specified for two model-selection type methods and one aggregation method. The approach has been applied to Gaussian, volatility, Poisson, exponential and Bernoulli models. A universal procedure for the choice of critical values has been developed. The influence of global parameters such as the test level on the critical values has been investigated. Theoretical properties of the procedures have been investigated in the most unified way possible so far. It is hardly possible to pick the best method among the considered ones. With a good choice of critical values all methods considered demonstrate very reasonable performance on simulated and real-life data.

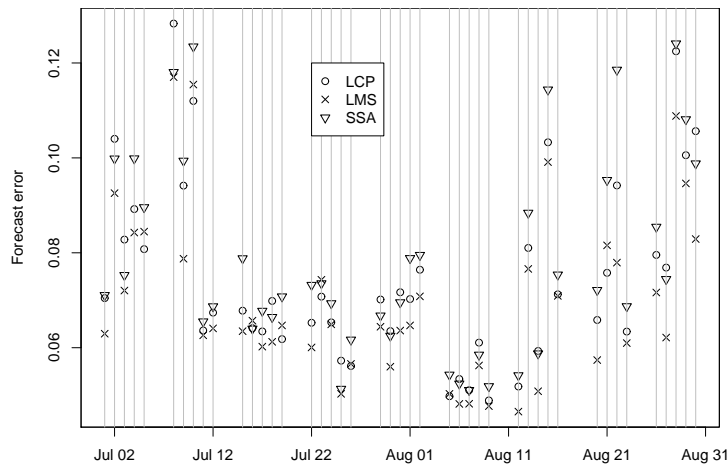

Fig 15. *Median forecasting error of the estimation of times between transactions involving the BHP Billiton stock on ASX.*

# Appendix

**Lemma 2 (Information-theoretical bound).** *Let $\mathbf{P}$ and $\mathbf{P}_0$ be two measures such that the Kullback – Leibler divergence $\mathbf{E}\log(d\mathbf{P}/d\mathbf{P}_0)$, satisfies*

$$\mathbf{E}\log(d\mathbf{P}/d\mathbf{P}_0) \le \Delta < \infty.$$

*Then for any random variable $\zeta$ with $\mathbf{E}_0\zeta < \infty$*

$$\mathbf{E}\log(1+\zeta) \le \Delta + \mathbf{E}_0\zeta.$$

**Proof.** One can check that for any fixed $y$ the maximum of the function $f(x) = xy - x\log x + x$ is attained at $x = e^y$. Substituting $x = e^y$ one obtains the inequality $xy \le x\log x - x + e^y$. Using this inequality and the representation $\mathbf{E}\log(1+\zeta) = \mathbf{E}_0\{Z\log(1+\zeta)\}$ with $Z = d\mathbf{P}/d\mathbf{P}_0$, one obtains

$$\begin{aligned}\mathbf{E}\log(1+\zeta) &= \mathbf{E}_0\{Z\log(1+\zeta)\} \\ &\le \mathbf{E}_0(Z\log Z - Z) + \mathbf{E}_0(1+\zeta) \\ &= \mathbf{E}_0(Z\log Z) + \mathbf{E}_0\zeta - \mathbf{E}_0 Z + 1.\end{aligned}$$

It remains to note that $\mathbf{E}_0 Z = 1$ and $\mathbf{E}_0(Z\log Z) = \mathbf{E}\log Z = \mathbf{E}\log(d\mathbf{P}/d\mathbf{P}_0) \le \Delta$. □

**Lemma 3.** *(Polzehl and Spokoiny, 2006, Lemma 5.2) For distributions from the exponential family with the functions $B(\theta)$ and $C(\theta)$ (see (3)) continuously differentiable on $\Theta$, where $\Theta$ satisfies Condition 4, it holds for every $\theta_0, \theta_1, \theta_2 \in \Theta$*

$$\sqrt{\mathcal{K}(\theta_1,\theta_2)} \le \mathfrak{a}\left(\sqrt{\mathcal{K}(\theta_1,\theta_0)} + \sqrt{\mathcal{K}(\theta_2,\theta_0)}\right).$$

*Moreover, for any sequence* $\theta_0, \theta_1, \ldots, \theta_m$

$$\sqrt{\mathcal{K}(\theta_0, \theta_m)} \leq \mathfrak{a}\left(\sqrt{\mathcal{K}(\theta_0, \theta_1)} + \cdots + \sqrt{\mathcal{K}(\theta_{m-1}, \theta_m)}\right).$$